\documentclass[preprint,3p,12pt]{elsarticle}
\biboptions{sort&compress}

\usepackage{amsmath}
\usepackage{graphicx}% Include figure files
\usepackage{dcolumn}% Align table columns on decimal point
\usepackage{bm}% bold math
\usepackage{caption}
\usepackage{subfigure}

\setcitestyle{super,open={},close={}}

\begin{document}

%\captionsetup{tablename=Table}
%\renewcommand{\thetable}{\arabic{table}}

%\preprint{APS/123-QED}
\begin{frontmatter}

\title{On the explicit representation of orthonormal Bernstein Polynomials}% Force line breaks with \\

\author{Michael A. Bellucci \corref{cor1}}
\ead{bellucci@mit.edu}

\cortext[cor1]{Corresponding author}
%\cortext[cor2]{Principal corresponding author}
%\fntext[fn1]{This is the specimen author footnote.}
%\fntext[fn2]{Another author footnote, but a little more longer.}
\address{Department of Chemical Engineering, Massachusetts Institute of Technology, Cambridge, Massachusetts 02139, U.S.A.}

\begin{abstract}
In this work we present an explicit representation of the orthonormal Bernstein polynomials and demonstrate that they can be generated from a linear combination of non-orthonormal Bernstein polynomials.  In addition, we report a set of $n$ Sturm-Liouville eigenvalue equations, where each of the $n$ eigenvalue equations have the orthonormal Bernstein polynomials of degree $n$ as their solution set.  We also show that each of the $n$ Sturm-Liouville operators are naturally self-adjoint.  While the orthonormal Bernstein polynomials can be used in a variety of different applications, we demonstrate the utility of these polynomials here by using them in a generalized Fourier series to approximate curves and surfaces.  Using the orthonormal Bernstein polynomial basis, we show that highly accurate approximations to curves and surfaces can be obtained by using small sized basis sets.  Finally, we demonstrate how the orthonormal Bernstein polynomials can be used to find the set of control points of B\'{e}zier curves or B\'{e}zier surfaces that best approximate a function. 

\end{abstract}


\end{frontmatter}

%\maketitle

\section{Introduction}
Bernstein polynomials are of great practical importance in the field of computer aided-geometric design as well as numerous other fields of mathematics because of their many useful properties.\cite{Farouki_2012, Farouki_1988,Bohm_1984,Farin_2002,Goldman_2003,Hoschek_1993,Prautzsch_2002}  Perhaps the best known practical use of Bernstein polynomials is in the definition of B\'{e}zier curves and B\'{e}zier surfaces, which are parametric curves and surfaces that use a Bernstein polynomial basis set in their representation.  B\'{e}zier curves and surfaces can be used to approximate any curve or surface to a high degree of accuracy, and therefore, they are very important tools used in computer graphics.\cite{Bohm_1984,Farin_2002,Goldman_2003,Hoschek_1993,Prautzsch_2002,Hormann_2008,Sederberg_1986}  However, Bernstein polynomials have numerous other applications aside from computer graphics.  Bernstein polynomials have been used in Galerkin methods and collocation methods to solve elliptic and hyperbolic partial differential equations.\cite{Bhatti_2007,Bhatta_2006,Doha_2011,Doha_2010,Mirkov_2013}  In addition, they are fundamental to approximation theory as they provide a way to prove the Weierstrass approximation theorem, which states that any continuous function on a closed and bounded interval can be uniformly approximated on that interval by polynomials to any degree of accuracy.\cite{Bernstein_1912,Weierstrass_1885}  Bernstein polynomials also have applications in optimal control theory,\cite{Yousefi_2010, Sanchooli_2010, Alipour_2013} stochastic dynamics,\cite{Kowalski_2006} and in the modeling of chemical reactions, where they can be used in B\'{e}zier curves to represent the most probable reaction path in high dimensional configuration space.\cite{Bellucci_2014}  

Despite the fact that Bernstein polynomials have many useful properties, one property they do not possess is orthogonality.  For many applications, such as least squares approximation and finite element methods, the orthogonality property is particularly useful, and as a result, the application of Bernstein polynomials in these methods is often less convenient than traditional orthogonal polynomials such as Legendre polynomials, Chebyshev polynomials, or Jacobi polynomials.  To overcome this difficulty, the Bernstein polynomial basis is often transformed into an orthogonal polynomial basis using a transformation matrix.\cite{Farouki_2012,Boyd_2008, Coluccio_2008, Rababah_2003, Rababah_2004, Rababah_2007, Farouki_2003, Doha_2010, Farouki_2000}  However, as the degree of the polynomial basis increases, the transformation matrix between basis sets can become ill-conditioned which can introduce substantial error into numerical calculations.\cite{Farouki_2012, Coluccio_2008, Farouki_2000, Hermann_1996}  Alternatively, the orthonormal Bernstein polynomial basis can be generated through a Gram-Schmidt orthonormalization process, but this process must be repeated every time the degree of the polynomial basis is increased.  It would clearly be beneficial in many of the applications discussed here to have an explicit representation to generate orthonormal Bernstein polynomials, but to the best of the authors' knowledge, there is no explicit representation of orthonormal Bernstein polynomials in the literature.  The aim of this paper is to present the explicit representation of the orthonormal Bernstein polynomials, discuss their corresponding Sturm-Liouville equation, and demonstrate their utility for curve and surface approximation.

\section{Bernstein Polynomials}
The Bernstein basis polynomials of degree $n$ form a complete basis over the interval $[0,1]$ and are defined by
\\
\begin{equation}
B_{j,n}(t)={n \choose j} t^j (1-t)^{n-j} \quad j=0,1,\ldots,n
\end{equation}
\\
where $t$ is a parameter.  However, the Bernstein basis polynomials can be generalized to cover an arbitrary interval $[a,b$] by normalizing $t$ over the interval $[a,b]$, \textit{i.e.} $t=(x-a)/(b-a)$, which leads to the following
\\
\begin{equation}
B_{j,n}(x)={n \choose j} \frac{(x-a)^j (b-x)^{n-j}}{(b-a)^n} \quad j=0,1,\ldots,n.
\label{eq:gen_bc}
\end{equation}
\\
These polynomials satisfy symmetry $B_{j,n}(x)=B_{n-j,n}(1-x)$, positivity $B_{j,n}(x) \ge 0$, and form a partition of unity $\sum_{j=0}^{n} B_{j,n}(x)=1$ on the defining interval $[a,b]$.  Moreover, they satisfy a number of other useful properties \cite{Farouki_2012} that we do not discuss in detail here.  

By taking a linear combination of Bernstein polynomials we can define a generalized parametric curve over the interval $[a,b]$, which is known as the B\'{e}zier curve,
 \\
\begin{equation}
f(x)=\sum_{j=0}^n B_{j,n}(x) P_j,
\label{eq:bc}
\end{equation}
\\
where $P$ is a set of coefficients, commonly referred to as control points.  An $n$th degree B\'{e}zier curve consists of $n+1$ Bernstein polynomials, which form a basis for the linear space $V_n$ consisting of all polynomials of degree $m$, where $m \le n$.  Using (\ref{eq:bc}), we can accurately represent a function $f(x)$ in the interval $[a,b]$ by finding the set of control points, $P$, that best approximate the function $f(x)$.  Similarly, a generalized surface in $R^3$ over the arbitrary interval $[a,b] \times [c,d]$ can be defined by the tensor product of Bernstein basis polynomials using the following
\\
\begin{equation}
 f(x,y)=\sum_{i=0}^n \sum_{j=0}^m B_{i,n}(x) B_{j,m}(y) P_{i,j},
 \label{eq:bcs}
 \end{equation}
 \\
where $P$ is a control point matrix and the surface is defined by $(n+1)(m+1)$ Bernstein basis polynomials.

%\section{Function Approximation}
%A common and important problem in many different fields involves approximating some target function, $f(x)$, or data using a simple functional form or spline function.  Spline functions are particularly suited for this task, but often is the case where it is desirable to approximate     

\section{Orthonormal Bernstein Polynomials}
The explicit representation of the orthonormal Bernstein polynomials, denoted by $\phi_{j,n}(t)$ here, was discovered by analyzing the resulting orthonormal polynomials after applying the Gram-Schmidt process on sets of Bernstein polynomials of varying degree $n$.  For example, for $n$=5, using the Gram-Schmidt process on $B_{j,5}(t)$, normalizing, and simplifying the resulting functions, we get the following set of orthonormal polynomials 
\\
\begin{equation}
\phi_{0,5}(t)=\sqrt{11}(1-t)^5  
\nonumber
\end{equation}
\begin{equation}
\phi_{1,5}(t)=3(1-t)^4(11t-1)
\nonumber
\end{equation}
\begin{equation}
\phi_{2,5}(t)=\sqrt{7}(1-t)^3(55t^2-20t+1)
\end{equation}
\begin{equation}
\phi_{3,5}(t)=\sqrt{5}(1-t)^2(165t^3-135t^2+27t-1)
\nonumber
\end{equation}
\begin{equation}
\phi_{4,5}(t)=\sqrt{3}(1-t)(330t^4-480t^3+216t^2-32t+1)
\nonumber
\end{equation}
\begin{equation}
\phi_{5,5}(t)=462t^5-1050t^4+840t^3-280t^2+35t-1.
\nonumber
\end{equation}
\\
We can see from these equations that the orthonormal Bernstein polynomials are, in general, a product of a factorable polynomial and a non-factorable polynomial.  For the factorable part of these polynomials, there exists a pattern of the form
\\
\begin{equation}
\bigg{(}\sqrt{2(n-j)+1}\bigg{)} (1-t)^{n-j} \quad j=0,1,\ldots,n.
\end{equation}
\\
While it is less clear that there is a pattern in the non-factorable part of these polynomials, the pattern can be determined by analyzing the binomial coefficients present in Pascal's triangle.  In doing this, we have determined the explicit representation for the orthonormal Bernstein polynomials to be
\\
\begin{equation}
\phi_{j,n}(t)=\bigg{(}\sqrt{2(n-j)+1}\bigg{)} (1-t)^{n-j} \sum_{k=0}^j (-1)^k {2n+1-k \choose j-k} {j \choose k} t^{j-k}.
\label{eq:bc0_on}
\end{equation}
\\
In addition, (\ref{eq:bc0_on}) can be written in a simpler form in terms of the original non-orthonormal Bernstein basis functions as
\\
\begin{equation}
\phi_{j,n}(t)=\sqrt{2(n-j)+1}\sum_{k=0}^j (-1)^k \frac{{2n +1 -k \choose j-k}{j \choose k}}{{n-k \choose j-k}} B_{j-k,n-k}(t),
\label{eq:bc_on_b}
\end{equation}
\\
which is a remarkably simple formula that can be used to generate orthonormal Bernstein polynomials on the interval $[0,1]$.  To confirm the orthonormal relation of these polynomials, we first multiply two of these polynomials together and integrate to get
\\
\begin{equation}
\int_0^1\phi_{i,n}(t)\phi_{j,n}(t)dt=
\label{eq:bc_on_1}
\end{equation}
\begin{equation}
\sqrt{[2(n-i)+1][2(n-j)+1]}\sum_{k=0}^i\sum_{l=0}^j (-1)^{k+l} \frac{{2n +1 -k \choose i-k}{i \choose k} {2n +1 -l \choose j-l}{j \choose l} }{{n-k \choose i-k} {n-l \choose j-l}} \int_0^1 B_{i-k,n-k}(t) B_{j-l,n-l}(t)dt.
\nonumber
\end{equation}
\\
If we use the general relation \cite{Juttler_1998}
\\
\begin{equation}
\int_0^1 B_{p,q}(t) B_{r,s}(t)dt=\frac{{q \choose p}{s \choose r}}{\big{[}q+s+1\big{]}{q+s \choose p+r}},
\label{eq:bern_ip}
\end{equation}
\\
with $p=i-k$, $r=j-l$, $q=n-k$, $s=n-l$, and plug this into (\ref{eq:bc_on_1}) and simplify, we get the following
\\
\begin{equation}
\int_0^1\phi_{i,n}(t)\phi_{j,n}(t)dt=
\label{eq:bc_on_rel}
\end{equation}
\begin{equation}
\sqrt{[2(n-i)+1][2(n-j)+1]}\sum_{k=0}^i\sum_{l=0}^j (-1)^{k+l} \frac{{2n +1 -k \choose i-k}{i \choose k} {2n +1 -l \choose j-l}{j \choose l} }{ \big{[}2n+1-(k+l)\big{]}{2n-(k+l) \choose (i+j)-(k+l)} }
\nonumber
\end{equation}
\\

Using different values for $i$ and $j$, it is easy to verify the orthonormal relation of the polynomials in (\ref{eq:bc_on_b}) using (\ref{eq:bc_on_rel}).  Although the explicit representation presented in (\ref{eq:bc_on_b}) is for $t$ in the interval $[0,1]$, the orthonormal Bernstein polynomials on the arbitrary interval $[a,b]$ can easily be obtained by letting $t=(x-a)/(b-a)$ in (\ref{eq:bc_on_b}).  Finally, since the orthonormal polynomials can be generated from a Gram-Schmidt process, the orthonormal polynomials necessarily satisfy the following relations over the interval $[a,b]$
\\
\begin{equation}
\int_a^b \phi_{i,n}(x) B_{j,n}(x) dx=
\begin{cases}
(b-a)\sqrt{2(n-i)+1} \sum_{k=0}^i (-1)^k \frac{{2n+1-k \choose i-k} {i \choose k} {n \choose j}}{\big{[}2n+1-k\big{]} {2n-k \choose (i+j)-k}},& j \ge i\\
0,  & j < i
\end{cases}
\label{eq:bern_on_bern_ip}
\end{equation}
\\
where we have used (\ref{eq:bern_ip}) with $p=i-k$, $r=j$, $q=n-k$, and $s=n$.

\begin{figure}[t!]
\centering
\includegraphics[height=18pc,width=24.0pc,angle=0]{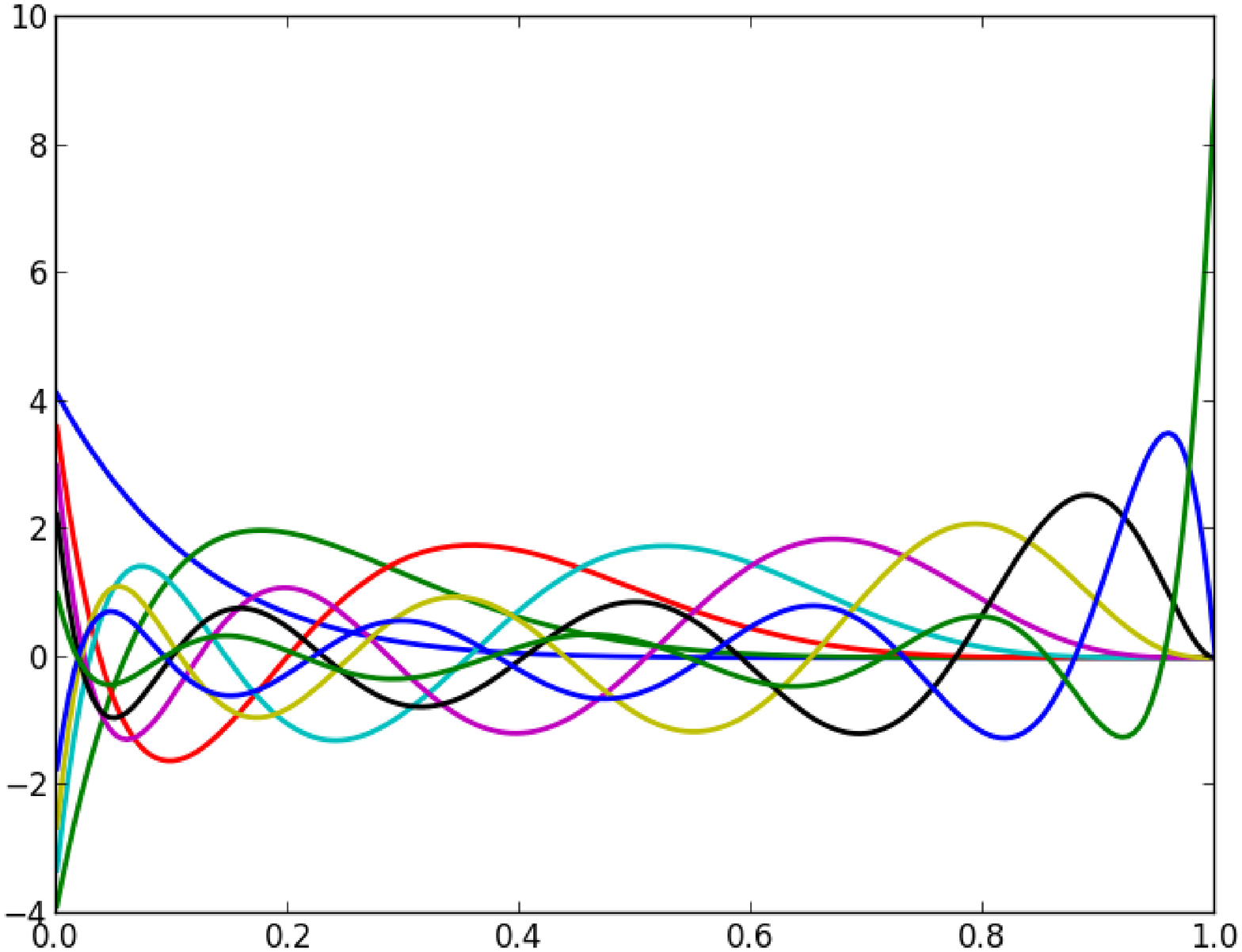}		
\caption{Orthonormal Bernstein polynomials with $n$=8}
\label{fig:on_bp}
\end{figure}

\section{Sturm-Liouville Equation}
The fact that the Bernstein polynomials in (\ref{eq:bc_on_b}) are orthonormal indicates that they can also be obtained as the solution of a Sturm-Liouville equation, since the solution set of Sturm-Liouville problems are orthogonal functions.  The Sturm-Liouville equation is a second-order linear differential equation of the form
\\
\begin{equation}
\frac{d}{dx}\bigg{[}p(x) \frac{d \phi}{dx}\bigg{]}+q(x)\phi+\lambda w(x) \phi=0,
\label{eq:SL_DE}
\end{equation}
\\
where $p(x)$, $q(x)$, and $w(x)$ are continuous and integrable real-valued functions on the finite interval $[a,b]$.  The solution set of this equation is a set of orthogonal functions, which are orthogonal with respect to the weight function $w(x)$, \textit{i.e.} they satisfy the following
\\
\begin{equation}
\int_a^b \phi_i \phi_j w(x) dx=\delta_{ij}
\end{equation}
\\
where $\delta_{ij}$ is the Kronecker delta function. The Sturm-Liouville equation can be simplified by defining a linear operator $\mathcal{L}$ on $C^2[a,b]$ as 
\\
\begin{equation}
\mathcal{L}\phi=\frac{d}{dx}\bigg{[}p(x) \frac{d \phi}{dx}\bigg{]}+q(x)\phi.
\end{equation}
\\
Using the operator $\mathcal{L}$, we can rewrite the Sturm-Liouville equation as an eigenvalue equation of the form
\\
\begin{equation}
\mathcal{L} \phi=-\lambda w(x) \phi.
\label{eq:SL_EP}
\end{equation}
\\
%There are many interesting properties of Sturm-Liouville systems, many of which we do not discuss here, but what is of particular importance to our discussion here is the self-adjoint property of the operator $\mathcal{L}$.  
Of central importance in the Sturm-Liouville theory is the self-adjoint property of the operator $\mathcal{L}$.  If the operator $\mathcal{L}$ is a self-adjoint operator with respect to the $L^2$ inner product space, then it can be shown that the resulting eigenfunctions, $\phi$, of $\mathcal{L}$ will be orthogonal and their corresponding eigenvalues, $\lambda$, will be real constants.  For the operator $\mathcal{L}$ to be self-adjoint, it must satisfy the following condition
\\
\begin{equation}
\big{<}u| \mathcal{L}v\big{>}-\big{<}\mathcal{L}u|v\big{>}=\int_a^b u^* \mathcal{L} v dx - \int_a^b (\mathcal{L} u)^* v dx=0,
\end{equation}
\\
where $u$ and $v$ are sufficiently smooth and integrable functions on the interval $[a,b]$.  Using integration by parts, it can be shown that the Sturm-Liouville operator is self-adjoint if and only if the following is satisfied
\\
\begin{equation}
\int_a^b u^* \mathcal{L} v dx - \int_a^b (\mathcal{L} u)^* v dx=p(x)\bigg{[}u^*\frac{dv}{dx}-v \frac{du^*}{dx}\bigg{]}\bigg{|}_a^b=0
\label{eq:SL_BC}
\end{equation}
\\
Therefore, in order for the eigenfunctions of (\ref{eq:SL_EP}) to be orthogonal and for the Sturm-Liouville eigenvalue problem to be well-posed, either the eigenfunctions must obey specific boundary conditions at the endpoints of the interval $[a,b]$, \textit{i.e.} homogenous boundary conditions, or the function $p(x)$ must vanish at the endpoints so that the right hand side of (\ref{eq:SL_BC}) is zero.  Often it is necessary to impose boundary conditions so that (\ref{eq:SL_BC}) is satisfied, but when $p(x)$ vanishes at the boundaries, it is not necessary to impose boundary conditions, aside from demanding that the eigenfunctions remain finite in the interval $[a,b]$.  This is case for the Sturm-Liouville eigenvalue equation of the Legendre polynomials, the Laguerre polynomials, and for the orthonormal Bernstein polynomials described here.  We have found that the orthonormal Bernstein polynomials, $\phi_{j,n}$, satisfy the following Sturm-Liouville eigenvalue equation
\\
\begin{equation}
% \frac{d}{dx}\bigg{[}x(1-x)^2 \frac{d \phi}{dx}\bigg{]}-n(n+2)x \phi+\lambda \phi=0,
 \frac{d}{dx}\bigg{[}x(1-x)^2 \frac{d \phi}{dx}\bigg{]}+n(n+2)(1-x) \phi+\lambda \phi=0,
\label{eq:SL_DE_1}
\end{equation}
\\
where $x$ is defined over the interval $[0,1]$ and the eigenvalues, $\lambda$, are defined by
\\
\begin{equation}
%\lambda=2\bigg{[}n\bigg{(}j+\frac{1}{2}\bigg{)}-{j \choose 2}\bigg{]} \quad j=0,1,\ldots,n.  
%\lambda=n(2j+1)-j(j-1)
\lambda=(n-j+1)(j-n)
\end{equation}
\\
Furthermore, since the sets of orthonormal Bernstein polynomials are distinctly different for each degree $n$, the operator $\mathcal{L}$ for (\ref{eq:SL_DE_1}) depends on the value $n$ through the function $q(x)=n(n+2)(1- x)$.  Therefore, (\ref{eq:SL_DE_1}) actually represents a set of $n$ Sturm-Liouville eigenvalue equations; one for each value of $n$.  For these equations, it is not necessary to impose boundary conditions since $p(x)=x(1-x)^2$ clearly vanishes at the endpoints of the interval $[0,1]$, demonstrating that the operator $\mathcal{L}$ associated with (\ref{eq:SL_DE_1}) is self-adjoint.

\section{Function Approximation}
While the orthonormal Bernstein polynomials can be used in many applications, we demonstrate their utility here in approximating curves and surfaces since function approximation is a common and important problem in many fields of applied mathematics and physics.  For a simple one dimensional curve, we can easily approximate the curve over the interval $[0,1]$ using the relation $g(t)=f(t)$, where $f(t)$ is the function to be approximated and $g(t)$ is a generalized Fourier series 
\\
\begin{equation}
 g(t)=\sum_{j=0}^n \phi_{j,n}(t) P_j.
\label{eq:bc_on}
\end{equation}
\\
This equation is analogous to the B\'{e}zier curve in (\ref{eq:bc}) but with an orthonormal Bernstein polynomials basis.  Approximating the curve amounts to choosing an appropriate value for the degree $n$ and finding the set of control points $[P_0,P_1,\ldots,P_n]$ that best fit the curve.  However, since the set of functions $[\phi_{0,n}(t),\phi_{1,n}(t),\ldots,\phi_{n,n}(t)]$ are orthonormal, the set of control points can easily be computed with
\\
\begin{equation}
P_j=\frac{\int_0^1 \phi_{j,n}(t) f(t) dt}{\int_0^1 \phi_{j,n}(t) \phi_{j,n}(t)dt},
\label{eq:bc_cp}
\end{equation}
\\
where the term in the denominator is equal to 1 due to the orthonormality of the basis functions.  In addition, if we let $t=(x-a)/(b-a)$ then we can approximate a curve over the arbitrary interval $[a,b]$.  In this case, the set of control points can be computed with
\\
\begin{equation}
P_j=\frac{1}{b-a}\int_a^b \phi_{j,n}(x) f(x) dx,
\label{eq:bc_on_P}
\end{equation}
\\
where the 1/(b-a) term comes from the following orthogonal relation
\\
\begin{equation}
\int_a^b \phi_{i,n}(x) \phi_{j,n}(x) dx=(b-a)\delta_{i,j}.
\end{equation}
\\
Similarly, to approximate a surface over the arbitrary interval $[a,b] \times [c,d]$, we can use a generalized Fourier series of the form
\\
\begin{equation}
 g(x,y)=\sum_{i=0}^n \sum_{j=0}^m \phi_{i,n}(x) \phi_{j,m}(y) P_{i,j},
 \label{eq:bcs_on}
\end{equation}
\\
and compute the elements of the control point matrix with the following
\\
\begin{equation}
P_{i,j}=\frac{1}{(b-a)(d-c)}\int_a^b \int_c^d f(x,y) \phi_{i,n}(x) \phi_{j,m}(y) dx dy.
\label{eq:bcs_on_P}
\end{equation}
\\
In addition, due to the relation in (\ref{eq:bern_on_bern_ip}), the orthonormal Bernstein polynomials can be used to find the control points of a B\'{e}zier curve or a B\'{e}zier surface of the forms presented in (\ref{eq:bc}) and (\ref{eq:bcs}), respectively.  This is immensely useful since it often can be difficult to determine the optimal control points that approximate a function using B\'{e}zier curves or surfaces, particularly for large degree $n$, due to the non-orthogonal properties of Bernstein polynomials.   For a B\'{e}zier curve, using the relation in (\ref{eq:bern_on_bern_ip}), we can find the control points over the arbitrary interval $[a,b]$ using a back substitution procedure with the following equation
\\
\begin{equation}
P_i=\frac{1}{\int_a^b \phi_{i,n}(x) B_{i,n}(x) dx}\bigg{(}\int_a^b f(x) \phi_{i,n}(x) dx- \sum_{j=i+1}^n P_j \int_a^b \phi_{i,n}(x) B_{j,n}(x) dx\bigg{)}
\label{eq:bc_P}
\end{equation}
\\
where the control points should be solved for in the order $i=n,n-1,\ldots,0$.  Moreover, this approach can be extended to solve for the control point matrix of a B\'{e}zier surface over the interval $[a,b]\times[c,d]$ with the following equation
\\
\begin{equation}
P_{i,j}=\frac{\int_a^b \int_c^d f(x,y) \phi_{i,n}(x) \phi_{j,m}(y) dx dy- \sum_{k=i}^n \sum_{l}^m P_{k,l} \int_a^b \phi_{i,n}(x) B_{k,n}(x) dx \int_c^d \phi_{j,m}(y) B_{l,m}(y) dy}{\int_a^b \phi_{i,n}(x) B_{i,n}(x) dx \int_c^d \phi_{j,m}(y) B_{j,m}(y) dy},
\nonumber
\end{equation}
\begin{equation}
\begin{cases}
l > j, & \text{if} \ k=i \\
l=j, & \text{otherwise}
\end{cases}
\label{eq:bcs_P}
\end{equation}
\\
where the control points should be solved for in the order $i=n,n-1, \ldots,0$ and $j=m,m-1,\ldots,0$.

To demonstrate the utility of these polynomials in approximating functions, we have chosen to approximate a parametric Lissajous curve, a sinc surface, and Langermann surface using the generalized Fourier series in (\ref{eq:bc_on}) and (\ref{eq:bcs_on}), respectively.  In addition, we have used the orthonormal Bernstein polynomials to find the best B\'{e}zier approximations to these test functions, \textit{i.e.} approximations that utilize non-orthonormal Bernstein polynomials. These test functions were chosen since they are common test functions used to challenge the capabilities of algorithms.  The Lissajous curve is defined by the parametric equations
\\
\begin{equation}
x(t)=A \sin(at+\delta) \quad y(t)=B \sin(bt).
\end{equation}
\\
The specific parameters chosen for this study was $A=1$, $B=1$, $a=4$, $b=3$, $\delta=\pi/3$ and $t$ was defined in the interval $[-\pi,\pi]$.  In addition, the sinc surface and Langermann surface we used were defined by the following equations
\\
\begin{equation}
f(x,y)=\frac{\sin(1.5\sqrt{x^2+y^2})}{s+1.5 \sqrt{x^2+y^2}}
\end{equation}
\begin{equation}
f(x,y)=\sum_{i=1}^p c_i \exp(-(x-q_i)^2/\pi-(y-r_i)^2/\pi)\cos(\pi (x-q_i)^2+\pi (y-r_i)^2)
\end{equation}
\\
where $s$ is a small constant to prevent dividing by zero at the origin, and $c_i$, $q_i$, and $r_i$ are parameters.  For the sinc surface, we chose $s=10^{-6}$ and the function was defined in the interval $[-8,8] \times [-8,8]$, whereas the Langermann surface was defined in the interval $[1,3] \times [1,3]$, and we chose $p=2$ with the following parameters 
\\
\begin{equation}
c=[1,2], \quad q=[2,3], \quad r=[3,2].
\end{equation}
\\

In order to find the best fit function, each of the test functions were discretized over a grid of points, $N$ points for the Lissajous curve and $N\times M$ points for the two surfaces.  Similarly, the orthonormal Bernstein polynomials and non-orthonormal Bernstein polynomials were discretized over these grids as well.  To find the control points for the B\'{e}zier approximations that utilized orthonormal Bernstein basis polynomials, the control points were computed with numerical integration using (\ref{eq:bc_on_P}) and (\ref{eq:bcs_on_P}) for the Lissajous curve and two surfaces, respectively.  For the B\'{e}zier approximations that utilized non-orthonormal Bernstein basis polynomials, the control points were computed with numerical integration using (\ref{eq:bc_P}) and (\ref{eq:bcs_P}).  To find the best approximation to the functions, in each case we gradually increased the degree of the polynomials until the error between the function and function approximation was a minimum.  The error for the curve and surfaces were measured with the following
\\
\begin{equation}
E=\frac{1}{N} \sum_{u=1}^N (x(t_u)-g(t_u))^2+(y(t_u)-h(t_u))^2
\end{equation}
\begin{equation}
E=\frac{1}{N M} \sum_{u=1}^N \sum_{v=1}^M (f(x_u,y_v)-g(x_u,y_v))^2.
\end{equation}
\\

\begin{figure}[t!]
    \centering
    \subfigure[]{\includegraphics[height=12.45pc,width=17.0pc,angle=0]{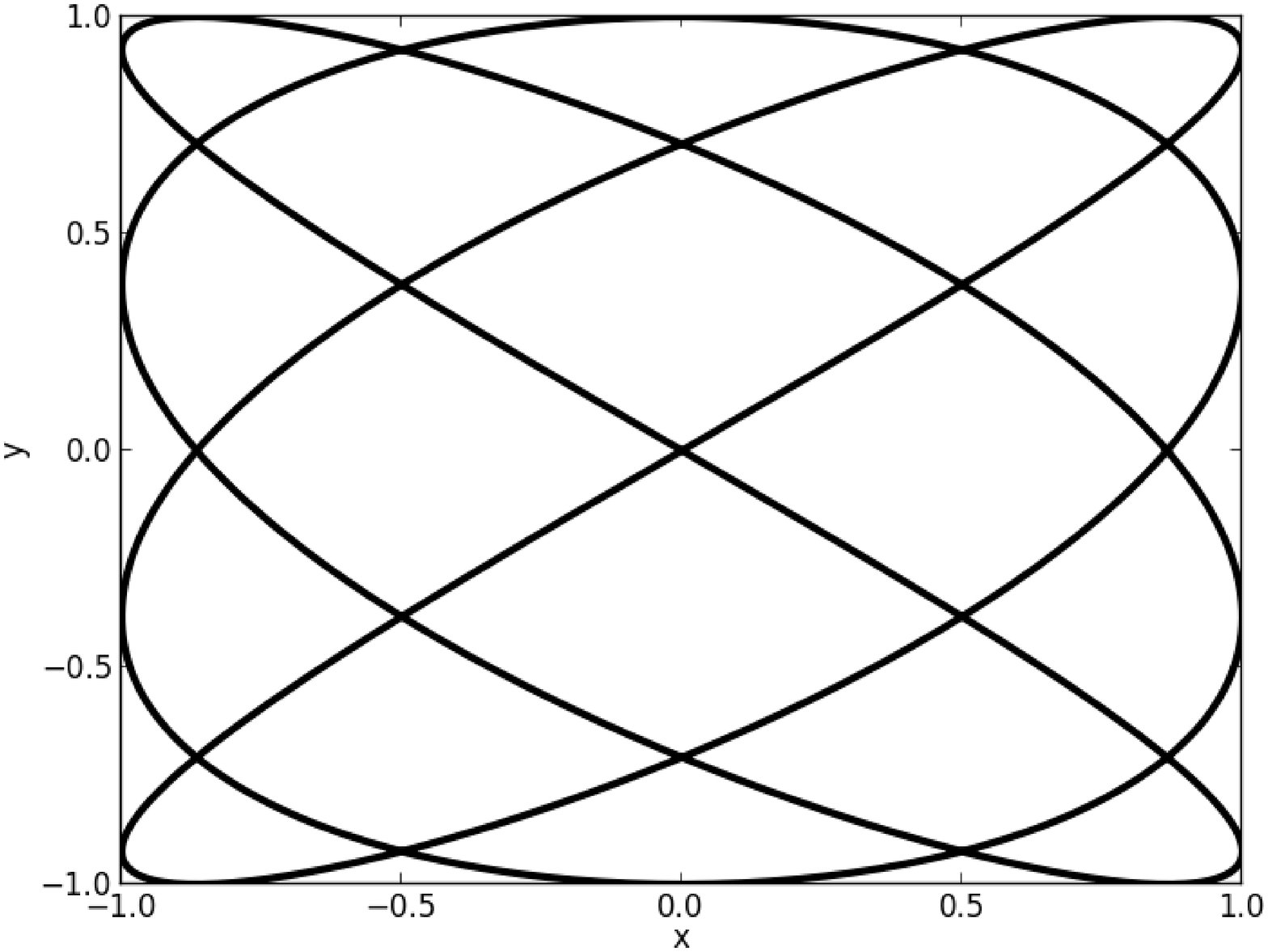}}		
    \subfigure[]{\includegraphics[height=12.45pc,width=17.0pc,angle=0]{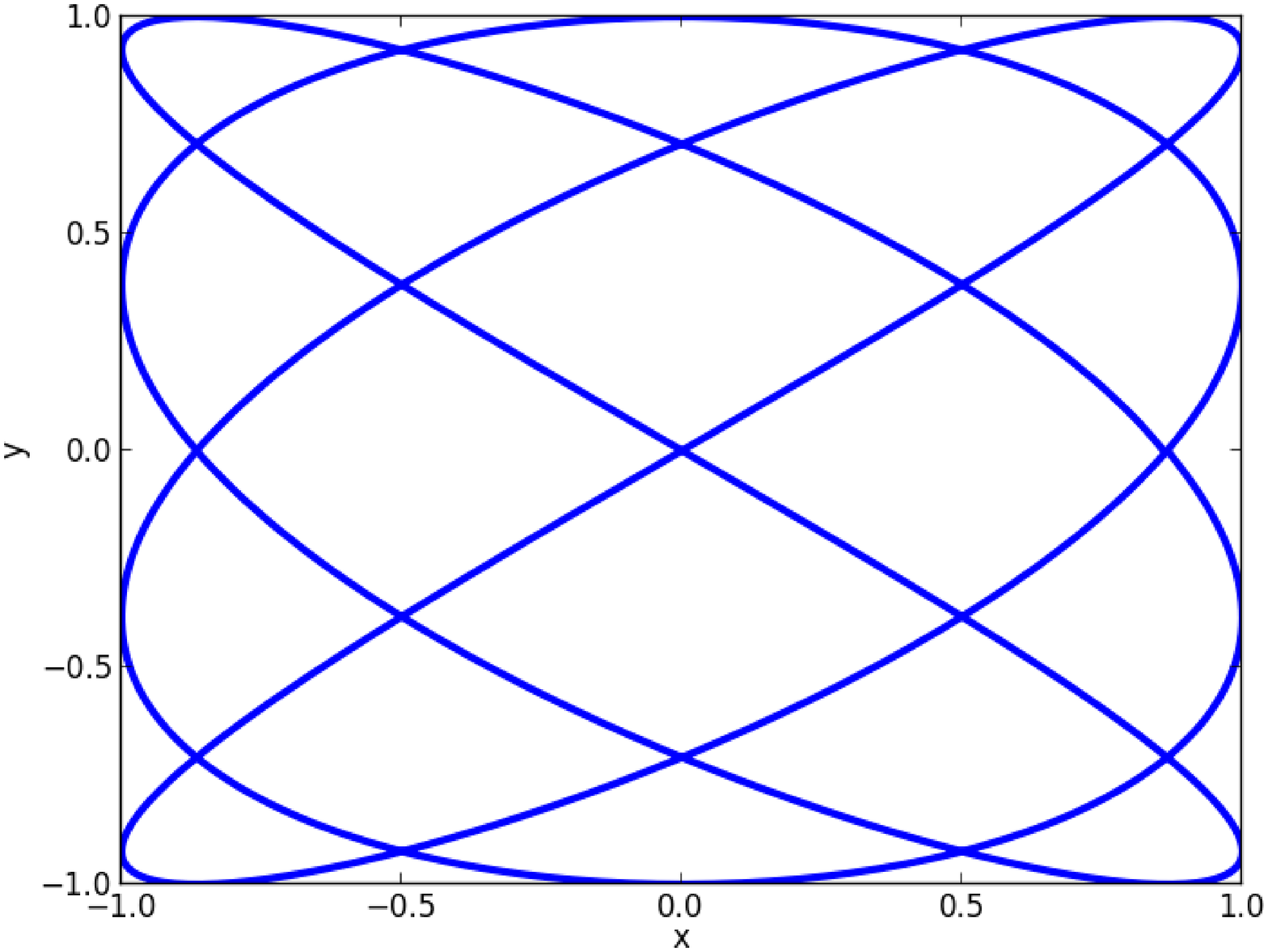}}		
    \subfigure[]{\includegraphics[height=12.45pc,width=17.0pc,angle=0]{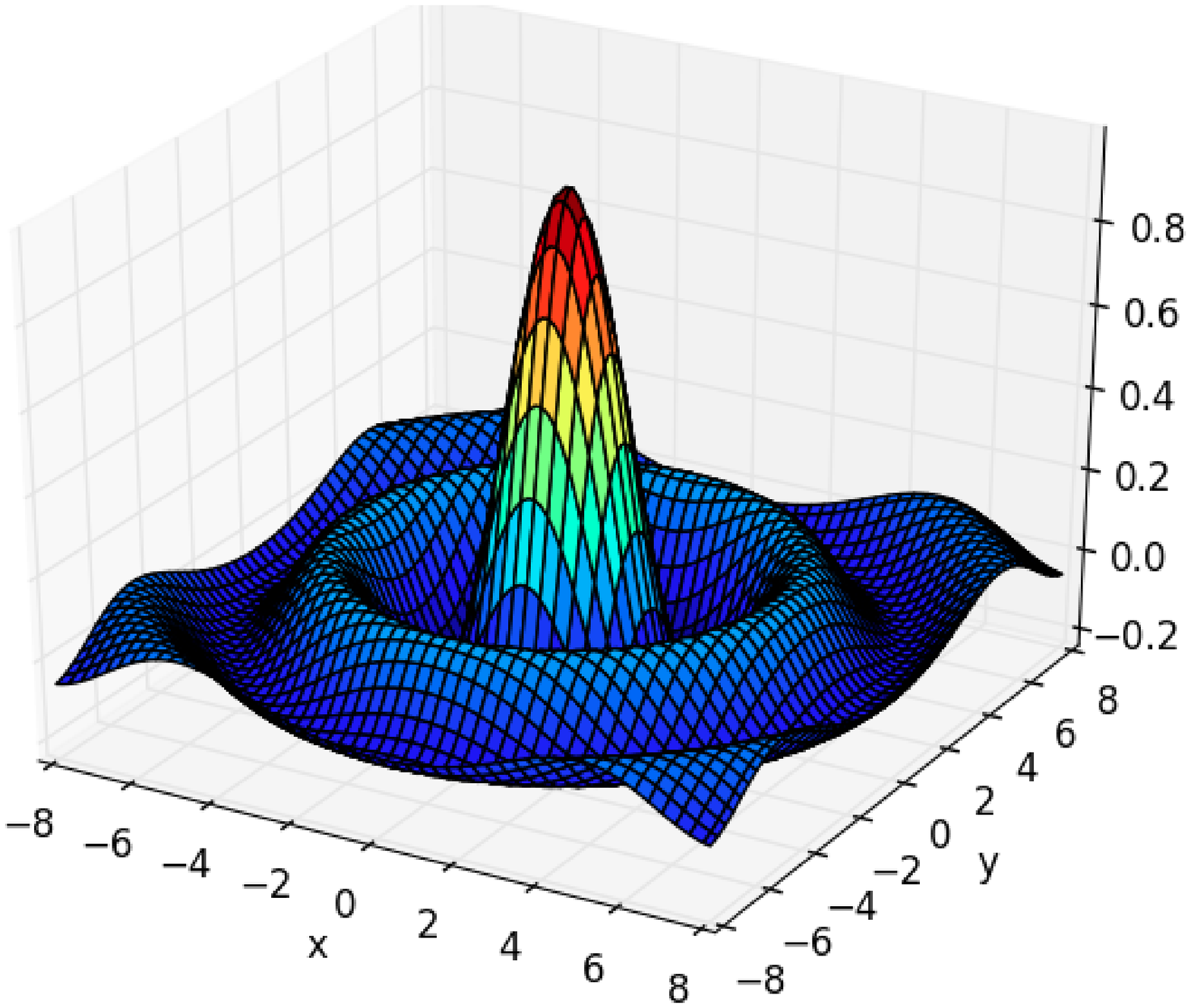}}		
    \subfigure[]{\includegraphics[height=12.45pc,width=17.0pc,angle=0]{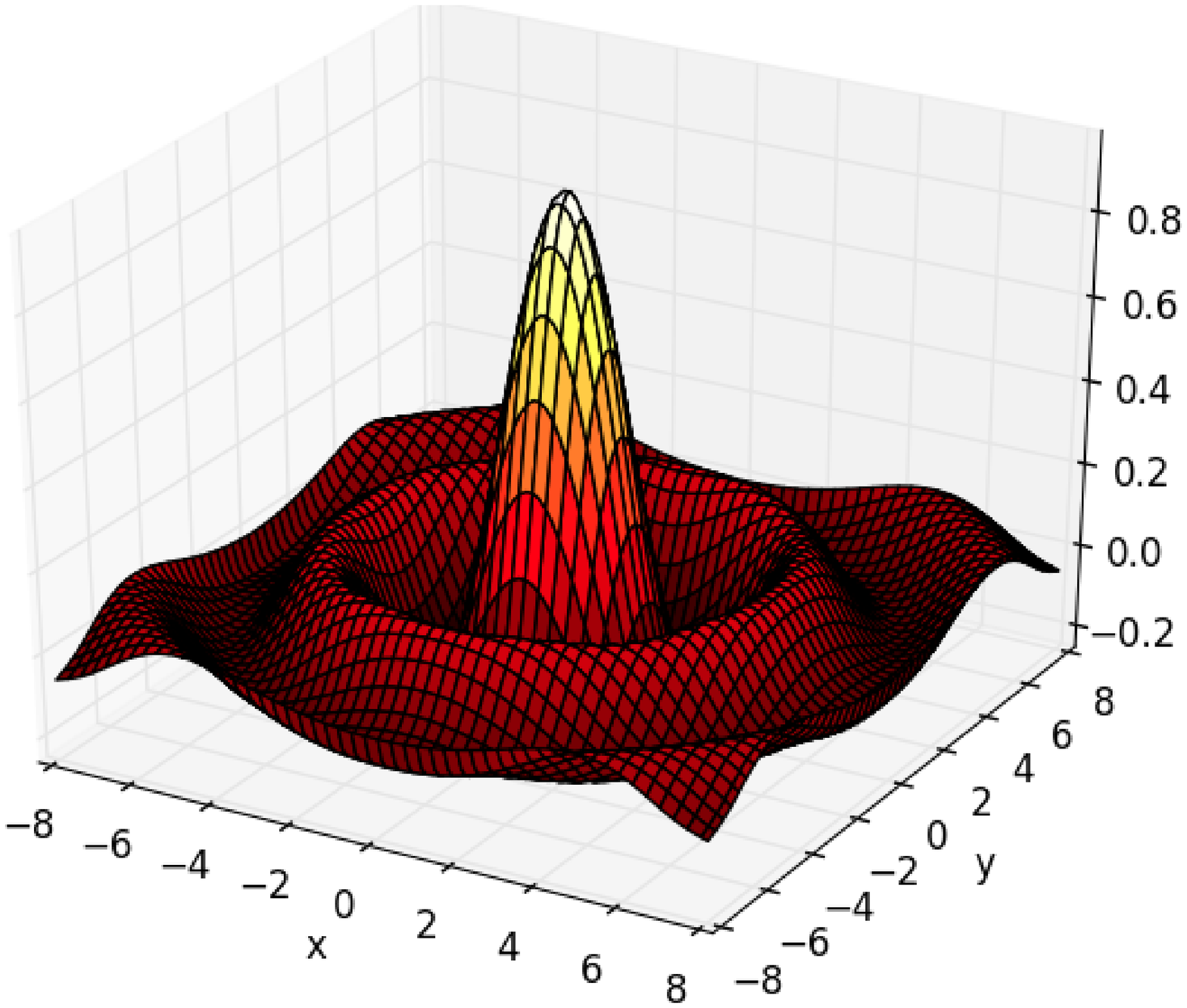}}		
    \subfigure[]{\includegraphics[height=12.45pc,width=17.0pc,angle=0]{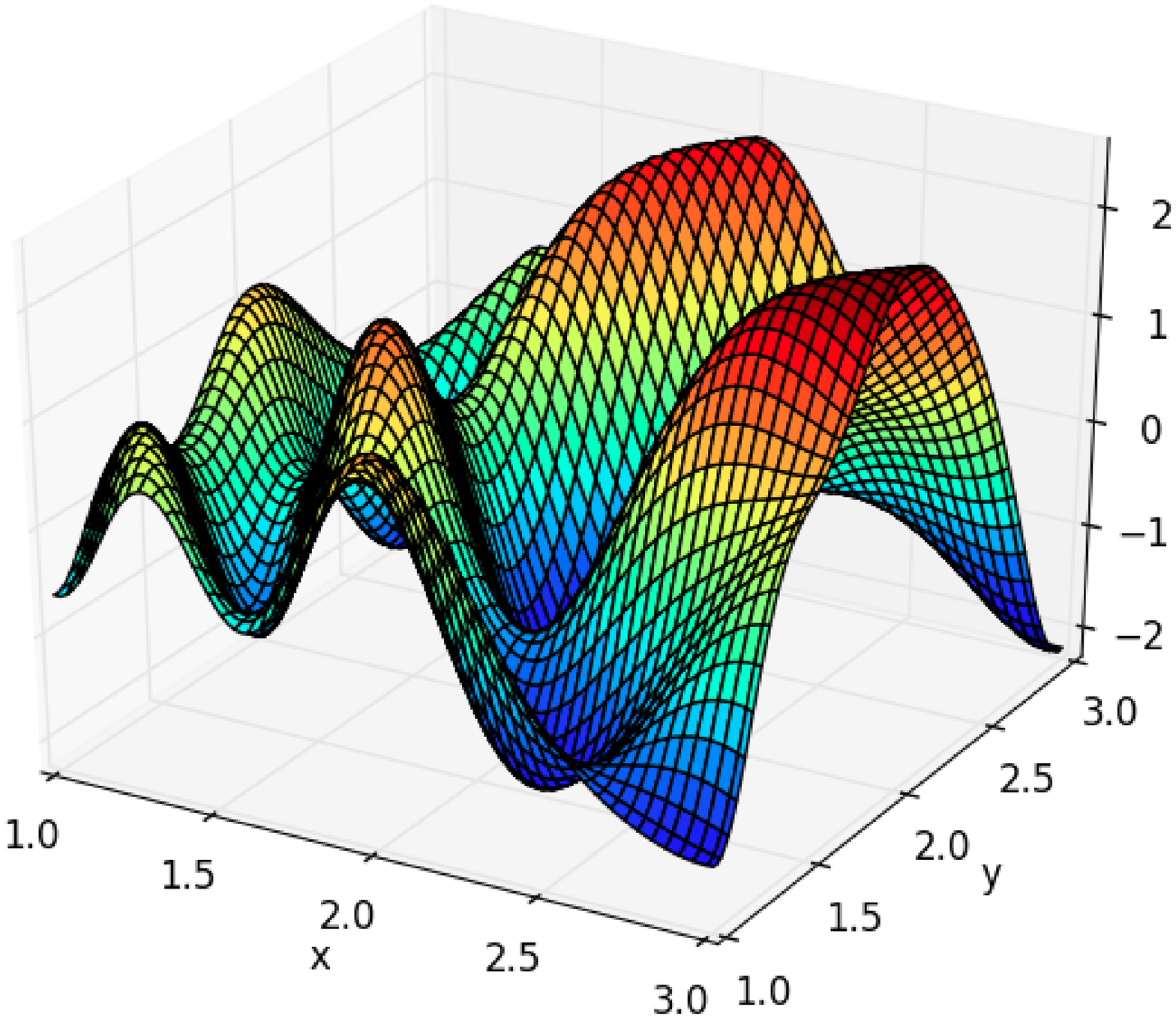}}		
    \subfigure[]{\includegraphics[height=12.45pc,width=17.0pc,angle=0]{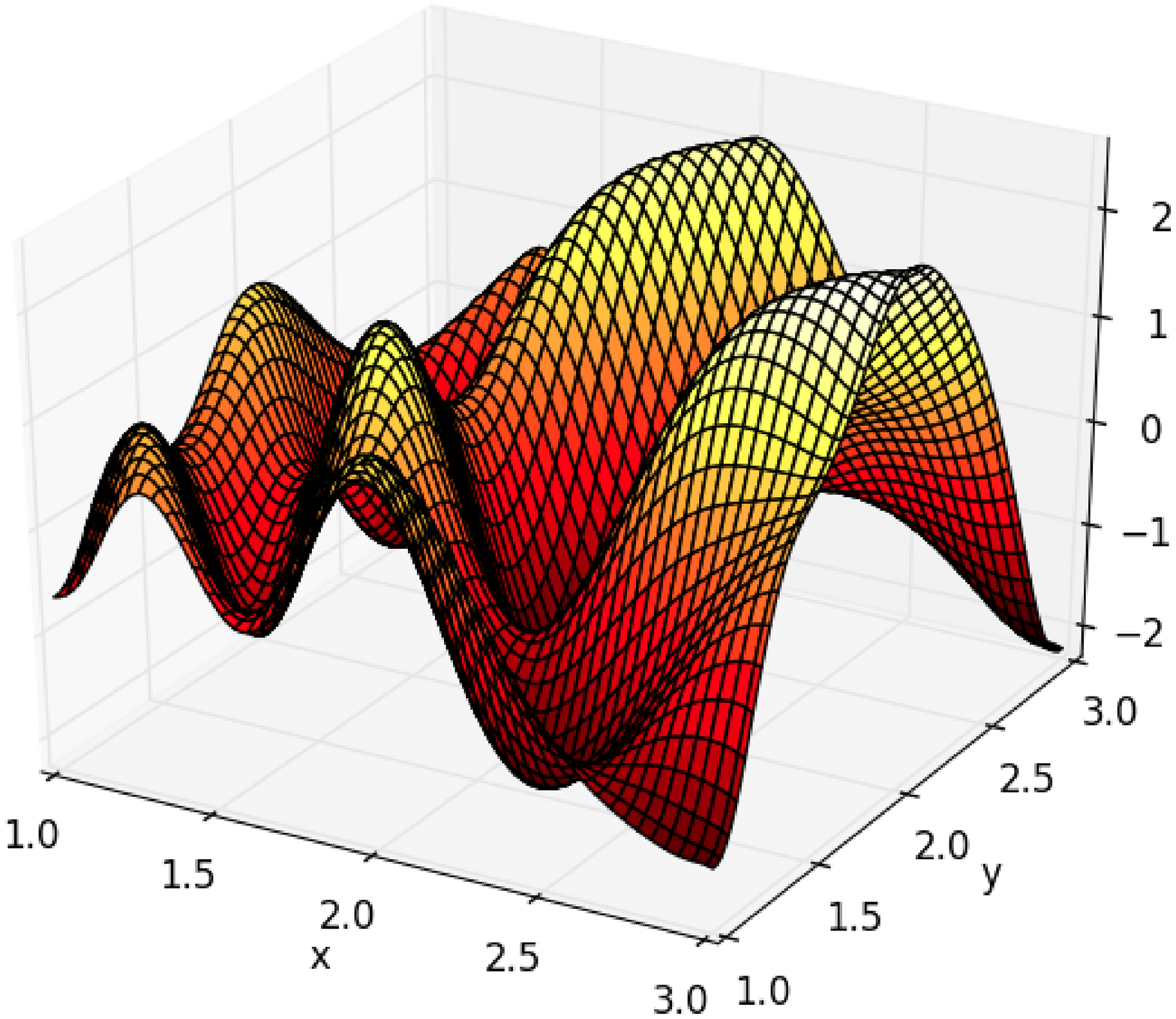}}		
    \caption{Upper panel: (a) Lissajous curve and (b) function approximation using degree $n=20$.  Mid Panel: (c) sinc surface and (d) function approximation using degree $n=12$ and $m=12$.  Lower Panel: (e) Langermann surface and (f) function approximation using degree $n=13$ and $m=13$. The function approximations shown here utilized orthonormal Bernstein polynomials.}
    \label{fig:fa}
\end{figure}

In Fig. \ref{fig:fa} the test functions and their B\'{e}zier curve and surface approximations are shown for comparison.  In this Figure, the approximations shown are for B\'{e}zier curve and surfaces that utilized orthonormal Bernstein polynomials.  We see that in all cases, the function approximations are very accurate representations of the original function.  For the Lissajous curve, the optimal function approximation was found with $n=20$,  which corresponded to an error of $2.1\times10^{-7}$ between the Lissajous curve and its B\'{e}zier curve approximation.  The large degree of the function approximation in this example has to do with the relatively high frequency of the sine waves in the Lissajous curve.  The orthonormal Bernstein polynomials resemble wave-like functions, and by increasing the degree of the polynomial, one effectively increases the frequency of these waves.  The optimal fit is found when the frequency of the basis functions are relatively close to the frequency of the sine waves in the Lissajous curve.  For lower frequency Lissajous curves, the optimal degree of the orthonormal Bernstein basis set can be significantly smaller.  For the B\'{e}zier curve approximation using non-orthonormal Bernstein polynomials, we found the optimal function approximation was once again found with $n=20$ and the error between the Lissajous curve and its B\'{e}zier curve approximation was $5.1 \times 10^{-8}$.  These results demonstrate that the optimal degree approximation is the same regardless of whether one uses orthonormal or non-orthonormal Bernstein polynomials as a basis, and we found this to be true for all test functions in this study.  In the mid-panel of Fig. \ref{fig:fa}, a comparison of the sinc surface and the optimal function approximation is shown.  For this surface, the best approximation was found with $n=12$ and $m=12$, which led to an error of $5.1\times 10^{-6}$.  Similarly, for the Langermann surface, the best approximation was found with $n=13$ and $m=13$, and the error between the surface and function approximation was $2.0 \times 10^{-6}$.  For surface approximations that utilized a non-orthonormal Bernstein basis set, the error between the surface approximation was $5.2 \times 10^{-6}$ and $2.3 \times 10^{-6}$ for the sinc surface and Langermann surface, respectively. These examples demonstrate that the orthonormal Bernstein polynomials are capable of reproducing complex functions with high accuracy using relatively low degree polynomial basis sets.  In addition, these results demonstrate that the orthonormal Bernstein polynomials can simplify the search for the optimal control points of B\'{e}zier approximations that use non-orthonormal Bernstein polynomials as a basis set.  The control points for the Lissajous curve and the two surfaces are given in the Supplementary Materials section.

\section{Conclusion}
In this work we have demonstrated that the orthonormal Bernstein polynomials can be generated from a linear combination of non-orthonormal Bernstein polynomials.  To the best of the authors' knowledge, this is the first explicit representation of the orthonormal Bernstein polynomials.  In addition, we have shown that the orthonormal Bernstein polynomials are the solution set of a set of $n$ Sturm-Liouville eigenvalue equations of the form in (\ref{eq:SL_DE}), where $p(x)=x(1-x)^2$, $q(x)=n(n+2)(1-x)$, $w(x)=1$, and $\lambda=(n-j+1)(j-n)$.  Moreover, we have demonstrated that the orthonormal Bernstein polynomials can be used in a generalized Fourier series to approximate curves and surfaces to a high degree of accuracy, and therefore, they can be very useful in computer-aided geometric design.  Furthermore, we have shown that they can simplify the search for the control points of B\'{e}zier curves and surfaces that best approximate functions.  However, these polynomials have many potential applications in numerous different fields of applied mathematics, where there is growing interest in the use of Bernstein polynomials in various applications.  We are currently applying these orthonormal polynomials in the modeling of chemical reactions to represent reaction paths in the high dimensional configuration space of chemical systems.  This will be the subject of our future work.  

\section*{Acknowledgements}
The author would like to thank Bernhardt Trout and Geoff Wood for useful discussions.  The author would like to kindly acknowledge support from Novartis through the Novartis-MIT Center for continuous manufacturing.

 %\section*{References}


\begin{thebibliography}{10}
\expandafter\ifx\csname url\endcsname\relax
  \def\url#1{\texttt{#1}}\fi
\expandafter\ifx\csname urlprefix\endcsname\relax\def\urlprefix{URL }\fi
\expandafter\ifx\csname href\endcsname\relax
  \def\href#1#2{#2} \def\path#1{#1}\fi

\bibitem{Farouki_2012}
R.~T. Farouki, The bernstein polynomial basis: A centennial retrospective,
  Comput. Aided Geom. D. 29~(6) (2012) 379--419.

\bibitem{Farouki_1988}
R.~T. Farouki, V.~T. Rajan, Algorithms for polynomials in bernstein form,
  Comput. Aided Geom. D. 5~(1) (1998) 1--26.

\bibitem{Bohm_1984}
W.~Bohm, G.~Farin, J.~Kahmann, A survey of curve and surface methods in cagd,
  Comput. Aided Geom. D. 1~(1) (1984) 1--60.

\bibitem{Farin_2002}
G.~E. Farin, Curves and Surfaces for CAGD: A Practical Guide (5th edition),
  Morgan Kaufmann, San Francisco, 2002.

\bibitem{Goldman_2003}
R.~Goldman, Pyramid Algorithms: A Dynamic Programming Approach to Curves and
  Surfaces for Geometric Modeling, Morgan Kaufmann, San Francisco, 2003.

\bibitem{Hoschek_1993}
J.~Hoschek, D.~Lasser, Fundamentals of Computer-Aided Geometric Design
  (translated by L. L. Schumaker), AK Peters, 1993.

\bibitem{Prautzsch_2002}
H.~Prautzsch, W.~Boehm, M.~Paluszny, B\'{e}zier and B-spline Techniques,
  Springer, Berlin, 2002.

\bibitem{Hormann_2008}
K.~Hormann, N.~Sukumar, Maximum entropy coordinates for arbitrary polytopes.,
  Computer Graphics Forum 27~(5) (2008) 1513--1520.

\bibitem{Sederberg_1986}
T.~W. Sederberg, S.~R. Parry, FreeÐform deformation of solid geometric models.,
  ACM SIGGRAPH Computer Graphics 20~(4) (1986) 151--160.

\bibitem{Bhatti_2007}
M.~I. Bhatti, P.~Bracken, Solutions of differential equations in a bernstein
  polynomial basis., J. Comput. Appl. Math 205~(1) (2012) 272--280.

\bibitem{Bhatta_2006}
D.~D. Bhatta, M.~I. Bhatti, Numerical solution of kdv equation using modified
  bernstein polynomials., Appl. Math. and Comput. 174~(2) (2006) 1255--1268.

\bibitem{Doha_2011}
E.~H. Doha, A.~H. Bhrawy, M.~A. Saker, Integrals of bernstein polynomials: An
  application for the solution of high even-order differential equations.,
  Appl. Math. Lett. 24~(1) (2011) 559--565.

\bibitem{Doha_2010}
E.~H. Doha, A.~H. Bhrawy, M.~A. Saker, On the derivatives of bernstein
  polynomials: An application for the solution of high even-order differential
  equations., Bound. Value Probl. 2011~(829543) (2011) 1--16.

\bibitem{Mirkov_2013}
N.~Mirkov, B.~Rasuo, Bernstein polynomial collocation method for elliptic
  boundary value problems., PAMM 13~(1) (2013) 421--422.

\bibitem{Bernstein_1912}
S.~N. Bernstein, D{\'e}monstration du th{\'e}or{\'e}me de weierstrass
  fond{\'e}e sur le calcul des probabilit{\'e}s., Communications de la
  Soci{\'e}t{\'e} Math{\'e}matique de Kharkov 2. Series XIII No. 1 (1912) 1--2.

\bibitem{Weierstrass_1885}
K.~Weierstrass, {\"U}ber die analytische Darstellbarkeit sogenannter
  willk{\"u}rlicher Functionen einer reellen Ver{\"a}nderlichen,
  Sitzungsberichte der K{\"o}niglich Preussischen Akademie der Wissenschaften
  zu Berlin, pp. 633-639 \& 789-805, reproduced in Mathematische Werke Vol. III,
  pp. 1-37, Georg Olms, Hildesheim., 1885.

\bibitem{Yousefi_2010}
S.~A. Yousefi, M.~Behroozifar, Operational matrices of bernstein polynomials
  and their applications., Int. J. Syst. Sci. 41~(6) (2010) 709--716.

\bibitem{Sanchooli_2010}
M.~Sanchooli, O.~S. Fard, Numerical scheme for fredholm integral equations
  optimal control problems via bernstein polynomials., Aust. J. Basic \& Appl.
  Sci. 4~(11) (2010) 5675--5682.

\bibitem{Alipour_2013}
M.~Alipour, D.~Rostamy, Bps operational matrices for solving time varying
  fractional optimal control problems., J. Math. Computer Sci 6 (2013)
  292--304.

\bibitem{Kowalski_2006}
E.~Kowalski, Bernstein polynomials and brownian motion., Am. Math. Mon.
  113~(10) (2006) 865--886.

\bibitem{Bellucci_2014}
M.~A. Bellucci, B.~L. Trout, B\'{e}zier curve string method for the study of
  rare events., In press.

\bibitem{Boyd_2008}
J.~P. Boyd, Exploiting parity in converting to and from bernstein polynomials
  and orthogonal polynomials., Appl. Math. Comput. 198~(2) (2008) 925--929.

\bibitem{Coluccio_2008}
L.~Coluccio, A.~Eisinberg, G.~Fedele, Gauss-lobatto to bernstein polynomials
  transformation., J. Comput. Appl. Math. 222~(2) (2008) 690--700.

\bibitem{Rababah_2003}
A.~Rababah, Transformation of chebyshevÐbernstein polynomial basis., Comput.
  Methods Appl. Math. 3~(4) (2003) 608--622.

\bibitem{Rababah_2004}
A.~Rababah, Jacobi-bernstein basis transformations., Comput. Methods Appl.
  Math. 4~(2) (2004) 206--214.

\bibitem{Rababah_2007}
A.~Rababah, M.~al~Natour, The weighted dual functionals for the univariate
  bernstein basis., Appl. Math. Comput. 186~(2) (2007) 1581--1590.

\bibitem{Farouki_2003}
R.~T. Farouki, T.~N.~T. Goodman, T.~Sauer, Construction of orthogonal bases for
  polynomials in bernstein form on triangular and simplex domains., Comput.
  Aided Geom. D. 20~(4) (2003) 209--230.

\bibitem{Farouki_2000}
R.~T. Farouki, LegendreÐbernstein basis transformations., J. Comput. Math.
  119~(1) (2000) 145--160.

\bibitem{Hermann_1996}
T.~Hermann, On the stability of polynomial transformations between taylor,
  bernstein, and hermite forms., Comput. Aided Geom. D. 13~(2) (1996) 307--320.

\bibitem{Juttler_1998}
B.~J{\"u}ttler, The dual basis functions for the bernstein polynomials, Adv.
  compute. Math. 8~(4) (1998) 345--352.

\end{thebibliography}
\end{document}